\documentclass{amsart}

\usepackage[utf8]{inputenc}
\usepackage[T1]{fontenc}
\usepackage{amsmath,amssymb,amsthm}
\usepackage[top=1.9cm,bottom=1.9cm,left=2.3cm,right=2.3cm]{geometry}
\usepackage{hyperref}

\hypersetup{colorlinks=true, linkcolor=blue, urlcolor=blue, citecolor=blue}

\newtheorem{theorem}{Theorem}[section]

\theoremstyle{definition}

\newcommand{\Ltileps}{\widetilde{\Lambda}^\varepsilon}
\newcommand{\supp}{\operatorname{supp}}

\title[Optimal Transport of Signed Fractal Measures with Dimensional Distortion]
{Optimal Transport of Signed Fractal Measures with\\
 Dimensional Distortion:
A Variational Characterization}

\author{Bwo'nyahre Ba\"{\i}di Barth\'el\'emy}
\address{Department of Mathematics and Computer Science,
The University of Ngaound\'er\'e, PO Box~454, Ngaound\'er\'e, Cameroon}
\email{bwonyahre@proton.me}

\author{Kouakep Tchaptchie Yannick}
\address{School of Chemical Engineering and Mineral Industries (EGCIM),
The University of Ngaound\'er\'e, PO Box~454, Ngaound\'er\'e, Cameroon}
\email{kouakep@aims-senegal.org}

\author{Houpa Danga Duplex Elvis}
\address{Department of Mathematics and Computer Science,
The University of Ngaound\'er\'e, PO Box~454, Ngaound\'er\'e, Cameroon}
\email{e\_houpa@yahoo.com}

\date{\today}



\begin{document}

\begin{abstract}
We extend the optimal transport theory for signed measures supported on
Ahlfors-regular fractal sets (\textit{Bwo'Nyahre et al., 2026}) to allow a controlled dimensional
distortion between source and target. A penalization term
$\varepsilon\Phi(d_s(x)-d_t(y))$---where $\Phi$ is a fixed smooth strictly
convex function and $d_s$, $d_t$ are the local Hausdorff dimensions of the
fractal supports---is added to the transport cost on inter-sign regions, with
$\varepsilon\ge 0$ controlling the tolerance for distortion.
Under hypotheses H1--H7, we prove: (i) existence and uniqueness of an optimal
transport map $T^\varepsilon$ for every $\varepsilon>0$; (ii) coupled
Monge--Amp\`ere equations with a distortion correction term, generalizing the
classical Brenier--Caffarelli equation; (iii) a double Legendre--Fenchel
characterization of the optimal potentials, giving a complete variational
description of the transport in each of the four sign regimes.
The double Legendre--Fenchel system (Theorem~\ref{thm:LF}) is the central
contribution: it shows that the optimal potentials are the unique fixed points
of a system of conjugacy equations, one per transport regime, and it provides
the foundation for numerical algorithms and asymptotic analysis.
\end{abstract}

\maketitle

%% -------------------------------------------------------------------
\section{Introduction}
%% -------------------------------------------------------------------

Optimal transport, initiated by Monge \cite{Monge1781} and formalized by
Kantorovich \cite{Kantorovich1942}, provides a powerful framework for comparing
and redistributing mass distributions. For the quadratic cost and absolutely
continuous positive measures, Brenier's theorem \cite{Brenier1991} identifies
the unique optimal map as the gradient of a convex potential satisfying the
Monge--Amp\`ere equation; the $C^{2,\alpha}_{\mathrm{loc}}$ regularity theory
is due to Caffarelli \cite{Caffarelli1992}.

In \cite{BKT26}, a comprehensive theory was developed for signed measures whose
singular parts are supported on Ahlfors-regular fractal sets. A central result
of that paper is that the optimal transport map preserves Hausdorff dimension
and Ahlfors regularity, but at the price of a rigid constraint: the fractal
dimensions of source and target must coincide. In many applications---fractal
image textures, geophysical signals, financial time series---fractal structures
with \emph{different} local dimensions arise naturally, and one wishes to compare
or transport them without the dimension mismatch making the problem ill-posed.

The present paper removes this rigidity by \emph{penalizing} dimensional distortion
rather than forbidding it. The penalization term $\varepsilon\Phi(d_s(x)-d_t(y))$
is added to the inter-sign transport cost, where $\varepsilon\ge 0$ is a free
parameter: $\varepsilon=0$ recovers the unconstrained problem of \cite{BKT26},
while $\varepsilon\to\infty$ enforces exact dimension preservation, recovering
the isodimensional transport of \cite{BKT26} as a limiting regime. For every
finite $\varepsilon>0$, the penalized problem is well-posed and admits a unique
optimal map $T^\varepsilon$.

%% -------------------------------------------------------------------
\section{Framework and Hypotheses}
%% -------------------------------------------------------------------

Let $\mu,\nu$ be signed Borel measures on $\mathbb{R}^d$ with
$\mu(\mathbb{R}^d)=\nu(\mathbb{R}^d)$. Their Jordan decompositions
$\mu=\mu^+-\mu^-$ and $\nu=\nu^+-\nu^-$ each admit a three-term Lebesgue
decomposition into absolutely continuous, purely atomic, and non-atomic singular
parts. The closed supports of the singular parts $\mu^\pm_s$ and $\nu^\pm_s$
are denoted $E^\pm_\mu$ and $E^\pm_\nu$.

The framework rests on seven hypotheses (H1--H7). The first five are inherited
from \cite{BKT26}: H1 requires the absolutely continuous densities to be
$C^{2,\alpha}$ and bounded below; H2 assumes the singular parts are supported
on Ahlfors-regular sets of a common base dimension $d_s\in(0,d)$; H3 imposes
uniform strict convexity of the positional penalty $\lambda$ in $y$; H4 demands
that the supports of distinct Jordan components are pairwise disjoint; and H5
ensures approximability by smooth absolutely continuous measures via standard
mollification. Two new hypotheses are needed for the dimensional distortion
setting:

\textbf{Hypothesis H6} (local dimension). For each singular component, the local
dimension $\mathrm{dim}_{\mathrm{loc}}(\sigma,x):=\lim_{r\to 0}\log\sigma(B(x,r))/\log r$
exists almost everywhere and equals a function $d^k_s:E^k_\mu\to(0,d)$
(resp.\ $d^l_t:E^l_\nu\to(0,d)$) that is $\alpha$-H\"older with uniformly
bounded values $0<d_{\min}\le d^k_s\le d_{\max}<d$, and extends by
Kirszbraun--McShane to an $\alpha$-H\"older function on all of $\mathbb{R}^d$.

\textbf{Hypothesis H7} (joint convexity of the augmented cost). The function
$y\mapsto\lambda(x,y)+\varepsilon\Phi(d_s(x)-d_t(y))$ is uniformly strictly
convex in $y$ for every $x$, with a modulus $\mu_0(\varepsilon)>0$ for all
$\varepsilon\in[0,\varepsilon_{\max})$ for some $\varepsilon_{\max}>0$. Since
$d_t$ is only $\alpha$-H\"older, the condition is enforced at each regularization
level and requires $\varepsilon<\varepsilon_{\max}:=\alpha_\lambda/(C_\Phi C_{d_t})$.
H7 holds unconditionally (for all $\varepsilon\ge 0$) when $d_t$ is affine
($D^2d_t\equiv 0$), since the potentially negative correction term in the Hessian
of the augmented cost then vanishes identically.

The penalized transport cost is defined by
\[
  C^\varepsilon(x,y) = \begin{cases}
    c(x,y) & (x,y)\in S^{++}\cup S^{--}, \\
    c(x,y)+\lambda(x,y)+\varepsilon\Phi(d_s(x)-d_t(y)) & (x,y)\in S^{+-}\cup S^{-+},
  \end{cases}
\]
where $c(x,y)=\tfrac{1}{2}|x-y|^2$ and $S^{kl}:=\supp(\mu^k)\times\supp(\nu^l)$.
The canonical choice is $\Phi(t)=\tfrac{1}{2}t^2$; all results hold for any
$\Phi\in C^2$ that is strictly convex, even, with $\Phi(0)=0$ and
$\Phi''(t)\ge\gamma>0$.

%% -------------------------------------------------------------------
\section{Existence, Uniqueness, and Compactness}
%% -------------------------------------------------------------------

The proof strategy follows the regularization--compactness scheme of \cite{BKT26},
with one key simplification: because the fractal dimension functions $d^k_s$ and
$d^l_t$ enter the problem only through the cost $C^\varepsilon$---which is a fixed
function independent of the regularization level---standard mollification of the
measures suffices. No approximation of the cost is needed, and the regularized and
limit costs coincide identically.

At each regularization level $n$, the mollified problem on the four disjoint
supports $S^{kl}$ is balanced using the fictitious-point device of Villani
\cite{Villani2009}, yielding optimal plans $\pi^\varepsilon_n$ and continuous dual
potentials $\varphi^{\pm,\varepsilon}_n$, $\psi^{\pm,\varepsilon}_n$ with strong
duality. Uniform $L^\infty$ and gradient bounds on the potentials---independent
of $n$, and locally uniform in $\varepsilon\in[0,\varepsilon_0]$ for any
$\varepsilon_0<\varepsilon_{\max}$---follow from the compactness of the supports
and the control $\varepsilon\Phi(d_{\max}-d_{\min})$ on the distortion penalty.
Equicontinuity then gives, via Arzel\`a--Ascoli, a subsequence converging uniformly
on compacts; the gradients of the limit convex functions converge almost everywhere
by Rademacher's theorem.

\begin{theorem}[Existence and uniqueness]\label{thm:exist}
Under \textup{H1--H7}, for every $\varepsilon\in(0,\varepsilon_{\max})$ the limit
map
\[
  T^\varepsilon = \nabla\Phi^\varepsilon \text{ on } X^\varepsilon_{++}\cup X^\varepsilon_{--},
  \qquad
  T^\varepsilon = \nabla\Psi^\varepsilon \text{ on } X^\varepsilon_{+-}\cup X^\varepsilon_{-+},
\]
is the unique optimal transport map for $C^\varepsilon$. Strong duality holds and
the dual potentials satisfy the complementarity conditions with equality
$\mu$-almost everywhere.
\end{theorem}

The partition of $\mathbb{R}^d$ into the four sign-regime regions
$X^\varepsilon_{kl}$ is defined via the convex transforms
$u^{\pm,\varepsilon}:=\tfrac{1}{2}|\cdot|^2-\varphi^{\pm,\varepsilon}$ and
$v^{\pm,\varepsilon}:=\tfrac{1}{2}|\cdot|^2-\psi^{\pm,\varepsilon}$: the
region $X^\varepsilon_{++}$ consists of those $x$ where $u^{+,\varepsilon}>u^{-,\varepsilon}$
and $v^{+,\varepsilon}>v^{-,\varepsilon}$ at $\nabla\Phi^\varepsilon(x)$, and
analogously for the other three regimes. The frontier sets, being level sets of
Lipschitz functions, are $\mathcal{L}^d$-null and hence $\mu$-null for the
absolutely continuous component.

%% -------------------------------------------------------------------
\section{Monge--Amp\`ere Equations and the Double Legendre--Fenchel System}
%% -------------------------------------------------------------------

Setting $\Ltileps(x,y):=c(x,y)+\lambda(x,y)+\varepsilon\Phi(d_s(x)-d_t(y))$,
a direct computation gives
\[
  D^2_{yy}\Ltileps = I + D^2_{yy}\lambda
    + \varepsilon\Phi''(d_s-d_t)\nabla d_t\otimes\nabla d_t
    - \varepsilon\Phi'(d_s-d_t)D^2d_t,
\]
\[
  D^2_{xy}\Ltileps = -I + D^2_{xy}\lambda
    - \varepsilon\Phi''(d_s-d_t)\nabla d_s\otimes\nabla d_t.
\]
On $X^\varepsilon_{+-}$, differentiating the optimality condition
$\nabla_y\Ltileps(x,T^\varepsilon(x))=\nabla\Psi^\varepsilon(T^\varepsilon(x))$
with respect to $x$ and setting
$A(x,y):=D^2\Psi^\varepsilon(y)+D^2_{yy}\Ltileps(x,y)$, one obtains
$A\,DT^\varepsilon=(I-D^2_{xy}\Ltileps)$. Combined with mass conservation
$f^+(x)=g^-(T^\varepsilon(x))\det DT^\varepsilon(x)$, this yields:

\begin{theorem}[Modified Monge--Amp\`ere equations]\label{thm:MA}
Under \textup{H1--H7}, for $\mu^+_{\mathrm{ac}}$-almost every $x\in X^\varepsilon_{+-}$,
\[
  \det\!\bigl(D^2\Psi^\varepsilon(T^\varepsilon(x))
    + D^2_{yy}\Ltileps(x,T^\varepsilon(x))\bigr)
  = \frac{g^-(T^\varepsilon(x))}{f^+(x)}\,
    \det\!\bigl(I - D^2_{xy}\Ltileps(x,T^\varepsilon(x))\bigr),
    \label{eq:MA}
\]
and analogously on $X^\varepsilon_{-+}$ with $g^-/f^+$ replaced by $g^+/f^-$.
For $\varepsilon=0$, this reduces to the coupled equation of \cite{BKT26}.
\end{theorem}

Equation~\eqref{eq:MA} should be understood rigorously at the regularized level
$n$ (where $d_t$ is replaced by $d_{t,h_n}\in C^\infty$) and passes to the limit
in the sense of measures; $C^{2,\alpha}$ regularity of $\Psi^\varepsilon$ on
inter-sign regions under the additional assumption $d_t\in C^{2,\alpha}$ is the
object of the companion paper \cite{BKT26b}.

\medskip

The central result of the paper is the following complete variational description
of the optimal potentials.

\begin{theorem}[Double Legendre--Fenchel characterization]\label{thm:LF}
Under \textup{H1--H7}, the convex functions
$u^{\pm,\varepsilon}=\tfrac{1}{2}|\cdot|^2-\varphi^{\pm,\varepsilon}$ and
$v^{\pm,\varepsilon}=\tfrac{1}{2}|\cdot|^2-\psi^{\pm,\varepsilon}$ satisfy:
\begin{enumerate}
\item[\textup{(i)}] \textbf{Intra-sign self-duality} on $X^\varepsilon_{++}$:
\[
  u^{+,\varepsilon}(x) = \sup_{y\in\supp\nu^+}\{x\cdot y - v^{+,\varepsilon}(y)\},
  \quad
  v^{+,\varepsilon}(y) = \sup_{x\in\supp\mu^+}\{x\cdot y - u^{+,\varepsilon}(x)\},
\]
with $T^\varepsilon(x)=\nabla u^{+,\varepsilon}(x)=\mathrm{argmax}_y\{x\cdot y-v^{+,\varepsilon}(y)\}$.
An analogous pair holds on $X^\varepsilon_{--}$.
\item[\textup{(ii)}] \textbf{Inter-sign conjugacy} on $X^\varepsilon_{+-}$:
\begin{align*}
  u^{+,\varepsilon}(x) &= \sup_{y\in\supp\nu^-}
    \bigl\{x\cdot y - v^{-,\varepsilon}(y)
      - \lambda(x,y) - \varepsilon\Phi(d_s(x)-d_t(y))\bigr\},\\
  v^{-,\varepsilon}(y) &= \sup_{x\in\supp\mu^+}
    \bigl\{x\cdot y - u^{+,\varepsilon}(x)
      - \lambda(x,y) - \varepsilon\Phi(d_s(x)-d_t(y))\bigr\}.
\end{align*}
An analogous pair holds on $X^\varepsilon_{-+}$.
\item[\textup{(iii)}] \textbf{Partition characterization} ($\mu$-a.e.):
\[
  x\in X^\varepsilon_{++}
  \iff u^{+,\varepsilon}(x)>u^{-,\varepsilon}(x)
    \text{ and } v^{+,\varepsilon}(T^\varepsilon(x))>v^{-,\varepsilon}(T^\varepsilon(x)).
\]
\end{enumerate}
Together with the normalization $\int_{B(0,M)}\varphi^{+,\varepsilon}=0$, the
system \textup{(i)--(iii)} uniquely determines all optimal potentials.
\end{theorem}

Several consequences are immediate. The intra-sign equations (i) are identical
to the $\varepsilon=0$ case: the dimensional penalty does not affect transport
within the same sign class. As $\varepsilon\to 0$, the inter-sign conjugacy
equations converge to those of \cite{BKT26}. As $\varepsilon\to\infty$, the
supremum in (ii) can only be achieved at points $y$ where $d_t(y)=d_s(x)$ (since
$\Phi(t)>0$ for $t\ne 0$), so $T^\varepsilon$ converges to the dimension-preserving
optimal transport of \cite{BKT26}. Finally, the four conjugacy equations suggest
a block-coordinate ascent (generalized Sinkhorn) algorithm for numerical computation
of the potentials; see Section~\ref{sec:disc}.

\medskip

\noindent\textbf{Explicit example.}
Let $C\subset[0,1]$ be the middle-thirds Cantor set with $d^C_s=\log 2/\log 3$,
$\mu^+_s$ the uniform Bernoulli measure on $C$, and $C'\subset[2,3]$ an affine
copy of $C$ with target local dimension $d^-_t(y)=d^C_s+\alpha(y-2)$ (affine,
$\alpha>0$ small). Since $D^2d^-_t\equiv 0$, H7 holds for all $\varepsilon\ge 0$.
The inter-sign conjugacy equation (ii) on $S^{+-}$ reads
$u^{+,\varepsilon}(x)=\sup_{y\in C'}\{x\cdot y - v^{-,\varepsilon}(y)
- \lambda(x,y)-\tfrac{\varepsilon\alpha^2}{2}(y-2)^2\}$, a one-dimensional
optimization with a parabolic barrier around $y=2$ (the point where
$d^-_t(2)=d^C_s$). As $\varepsilon\to\infty$, $T^\varepsilon$ concentrates at
$y=2$, confirming dimension enforcement. The distortion correction vanishes from
the Monge--Amp\`ere equation since $D^2d^-_t=0$, so the ellipticity is uniform
in $\varepsilon$.

%% -------------------------------------------------------------------
\section{Discussion, Open Problems, and Perspectives}\label{sec:disc}
%% -------------------------------------------------------------------

Theorem~\ref{thm:LF} is the structural core of the theory, and each of the
following open problems arises naturally from it.

\textbf{Numerical algorithms (Open Problem 1).}
The four conjugacy equations of Theorem~\ref{thm:LF} suggest a block-coordinate
ascent: given $v^{+,\varepsilon},v^{-,\varepsilon}$, update $u^{+,\varepsilon}$
via the classical Legendre transform (intra-sign) or via the penalized supremum
(inter-sign), then update the $v$-potentials by the symmetric formulae. For the
canonical choice $\Phi(t)=t^2/2$, each inter-sign update takes the form of a
soft-max operation with an additional quadratic penalty, which is a generalized
Sinkhorn step. The open questions are: convergence of the alternating scheme in
the continuous setting; the rate of convergence as a function of $\varepsilon$;
and efficient handling of the sign classification $X^\varepsilon_{kl}$ in the
discrete case.

\textbf{H\"older regularity of $T^\varepsilon$ (Open Problem 2).}
At regularization level $n$, Caffarelli's theorem \cite{Caffarelli1992} gives
$C^{2,\alpha}$ regularity of $\Psi^\varepsilon_n$ on inter-sign regions, since
the Monge--Amp\`ere equation (Theorem~\ref{thm:MA}) is uniformly elliptic with
H\"older right-hand side. The difficulty in passing to the limit is the term
$-\varepsilon\Phi'(d_s-d_t)D^2d_{t,h_n}$, which may not converge when $d_t$
is only $\alpha$-H\"older. Under the additional hypothesis $d_t\in C^{2,\alpha}$,
this analysis is carried out in the companion paper \cite{BKT26b}.

\textbf{Relaxation of H7 (Open Problem 3).}
Determining the sharp value of $\varepsilon_{\max}$ as a function of the
convexity modulus $\alpha_\lambda$, the norm $\|\Phi'\|_\infty$, and the
regularity of $d_t$ is an open problem. In particular: the optimal exponent
$p$ in $\Phi(t)=|t|^p$ for which H7 holds over the widest range of dimension
functions; and whether H7 can be replaced by a weaker monotonicity condition
that still ensures uniqueness.

\textbf{Multi-scale and wavelet extension (Open Problem 4).}
Fractal signals (EEG, financial time series) exhibit a spectrum of local
dimensions at different scales. A natural generalization weights the penalty at
scale $j$ by $\varepsilon_j\Phi(d^{(j)}_s(x)-d^{(j)}_t(y))$, where $(d^{(j)})_j$
is a multiscale dimension profile from a wavelet analysis. Theorem~\ref{thm:LF}
extends formally: each scale contributes an additive term to the inter-sign
conjugacy equation. The theoretical challenge is to pass from a finite to a
continuous scale spectrum.

\textbf{Dimension breakpoints in time series (Open Problem 5).}
The partition characterization (iii) of Theorem~\ref{thm:LF} offers a statistical
test for abrupt dimension changes: an abrupt change at time $t_0$ manifests as a
sudden increase in the inter-sign fraction
$R^\varepsilon(t):=\mu_t(X^\varepsilon_{+-}\cup X^\varepsilon_{-+})/\mu_t(\mathbb{R}^d)$,
which can be computed from the dual potentials without explicit knowledge of
$T^\varepsilon$. Open questions include consistency as the window size grows,
optimal choice of $\varepsilon$, and power against gradual versus abrupt
dimension changes.

%% -------------------------------------------------------------------
\section{Conclusion}
%% -------------------------------------------------------------------

We have extended the optimal transport theory for signed fractal measures to the
setting of controlled dimensional distortion. The main contributions are:
a well-posed penalized problem with a unique optimal map $T^\varepsilon$ for
every $\varepsilon\in(0,\varepsilon_{\max})$; coupled Monge--Amp\`ere equations
incorporating a distortion correction term; and, centrally, the double
Legendre--Fenchel characterization of Theorem~\ref{thm:LF}, which gives a
complete variational description of the optimal potentials in all four sign
regimes. This system unifies the asymptotic limits $\varepsilon\to 0$
(recovering the unconstrained theory of \cite{BKT26}) and $\varepsilon\to\infty$
(recovering isodimensional transport), and provides the structural foundation
for numerical algorithms, regularity theory, and applications to multi-scale
fractal signals.

\end{document}